\documentclass[mathematics,article,moreauthors,pdftex,10pt,a4paper]{mdpi} 
\firstpage{1} 
\doinum{10.3390/------}
\pubvolume{xx}
\pubyear{2017}
\copyrightyear{2017}
\history{Received: date; Accepted: date; Published: date}

\Title{A note on Hadamard fractional differential equations with varying coefficients and their applications in probability}

\Author{Roberto Garra $^{1 \ddagger,*}$, Enzo Orsingher $^{1 \ddagger}$ and Federico Polito $^{2 \ddagger}$}

\AuthorNames{Roberto Garra, Enzo Orsingher and Federico Polito}

\address{$^{1}$ \quad Dipartimento di Scienze Statistiche, ``Sapienza'' Universit\`{a} di Roma, P.~le A.\ Moro 5, 00185 Roma, Italy; \\
$^{2}$ \quad Dipartimento di Matematica ``G.~Peano'', Universit\`{a} degli Studi di Torino, Via Carlo Alberto 10, 10123 Torino, Italy.}

\corres{Correspondence: roberto.garra@sbai.uniroma1.it}

\abstract{In this paper we show several connections between special functions arising from generalized
	COM-Poisson-type statistical distributions and  integro-differential equations with varying coefficients involving Hadamard-type operators.
	New analytical results are obtained, showing the particular role of Hadamard-type derivatives in connection with a recently introduced
	generalization of the Le Roy function. 
	We are also able to prove a general connection between fractional hyper-Bessel-type equations involving Hadamard operators and Le Roy functions.}



\begin{document}
%%%%%%%%%%%%%%%%%%%%%%%%%%%%%%%%%%%%%%%%%%
%% Only for the journal Gels: Please place the Experimental Section after the Conclusions

%%%%%%%%%%%%%%%%%%%%%%%%%%%%%%%%%%%%%%%%%%
\setcounter{section}{0} %% Remove this when starting to work on the template.

\section{Introduction}

	The analysis of fractional differential equations involving Hadamard fractional derivatives has gained interest 
	in mathematical analysis, as proved, for example, by the publication of the recent monograph \cite{mono}.
    On the other hand, few results regarding the applications of Hadamard fractional differential equations in mathematical physics
    (see for example \cite{giusti}) and probability (see \cite{saxena}) exist.
    In \cite{had}, some analytical results regarding Hadamard fractional equations with time-varying coefficients have been pointed out.
    In particular, the following $\alpha$-Mittag--Leffler function was introduced:
    \begin{equation}\label{0}
    	E_{\alpha; \nu, \gamma}(z) =\sum_{k=0}^\infty \frac{z^k}{\left[\Gamma(\nu k+\gamma)\right]^\alpha},
    	\quad z\in \mathbb{C}, \alpha, \nu, \gamma \in \mathbb{C}.
    \end{equation}
    The $\alpha$-Mittag--Leffler function is an entire function of the complex variable $z$ if the parameters are such that $\Re (\nu) > 0 $,
    $\gamma \in \mathbb{R}$ and $\alpha\in \mathbb{R}^+$ (see \cite{main}). This function was independently introduced and studied by
    Gerhold in \cite{ge}. 
    In \cite{had}, we called this special function $\alpha$-Mittag--Leffler function, since it includes for $\alpha = 1$ the well-known
    two-parameters Mittag--Leffler function
    \begin{equation}
        E_{\nu, \gamma}(z) =\sum_{k=0}^\infty \frac{z^k}{\Gamma(\nu k+\gamma)}, \quad z\in \mathbb{C}, \ \nu \in \mathbb{C},
        \ \Re{(\nu)}>0, \ \gamma\in \mathbb{R}, 
    \end{equation}
    widely used in the theory of fractional differential equations (see the recent monograph \cite{gorenflo} and the references therein).
    Moreover, the $\alpha$-Mittag--Leffler function is a generalization of the so-called Le Roy function \cite{leroy}
    \begin{equation}
    	R_{\rho}(z)= \sum_{k=0}^\infty \frac{z^k}{[(k+1)!]^\alpha}, \quad z \in \mathbb{C}, \ \alpha >0.
    \end{equation}  
	In the more recent paper \cite{main}, the authors studied the asymptotic behavior and numerical simulation of this new class of special functions.
    
    We observe that the Le Roy functions are used in probability in the context of the studies of
    COM-Poisson distributions \cite{con}, which are special classes of weighted Poisson distributions (see for example \cite{Bala}).
    Therefore, the $\alpha$-Mittag--Leffler functions can also be used in the construction of a new generalization of the COM-Poisson
    distribution that can be interesting for statistical applications and in physics in the context of generalized coherent states \cite{Person}.
      
    The aim of this paper is to study some particular classes of Hadamard fractional integrals or differential equations whose
    solutions can be written in terms of the $\alpha$-Mittag--Leffler function \eqref{0} and somehow related
    to a \textit{fractional}-type generalization of the COM-Poisson distribution.           		
               		
    We also observe that Imoto \cite{imoto} has recently introduced the following generalization of the COM-Poisson distribution,
    \begin{equation}\label{imoo}
    	P\{N(t)= k\} = \frac{\left(\Gamma(\nu+k)\right)^r t^k}{k!C(r,\nu, t)}, \quad k\geq 0,
    \end{equation}
    involving the normalizing function 
	\begin{equation}\label{imo}
    	C(r,\nu, t)= \sum_{k=0}^\infty \frac{\left(\Gamma(\nu+k)\right)^r t^k}{k!},
    \end{equation}
    with $r<1/2$, $t>0$ and $1>\nu>0$ or $r=1$, $\nu = 1$ and $|t|<1$. The special function \eqref{imo} is somehow related to the
    generalization of the Le Roy function, while the distribution \eqref{imoo} includes the COM-Poisson for
    $\nu = 1$ and $r = 1-n$, $n\in \mathbb{R}$. In Section 4, we show that this function is related to integral equations with a time-varying
    coefficient involving Hadamard integrals.
   
	In Section 5, we present other results concerning the relation between Le Roy-type functions and Hadamard fractional differential equations.
   	We introduce a wide class of integro-differential equations extending the hyper-Bessel equations. This is a new interesting approach
   	in the context of the mathematical studies of fractional Bessel equations (we refer for example to the recent paper \cite{sitnik}
   	and the references therein).
   			     
    Concluding, the main aim of this paper is to establish a connection between some generalizations of the COM-Poisson distributions
    and integro-differential equations with time-varying coefficients involving Hadamard integrals or derivatives.
    As a by-product we suggest a possible application of the $\alpha$-Mittag--Leffler function to build a generalized COM-Poisson distribution that
    in future should be investigated in more detail.
 
%%%%%%%%%%%%%%%%%%%%%%%%%%%%%%%%%%%%%%%%%%

\section{Preliminaries about fractional Hadamard derivatives and integrals}

	Starting from the seminal paper by Hadamard \cite{Hadamard}, many papers have been devoted to the analysis of fractional operators with
	logarithmic kernels (we refer in particular to \cite{kilh}). In this section we briefly recall the definitions and main properties of
	Hadamard fractional integrals and derivatives and their Caputo-like regularizations recently introduced in the literature.
	\begin{Definition}
		\label{Integral}
		Let $t\in \mathbb{R}^+$ and $\Re (\alpha) >0$. The Hadamard fractional integral of order $\alpha$,
		applied to the function $f \in L^p[a,b]$, $1\leq p<+\infty$, $0<a<b<\infty$, for $t \in [a,b]$,
		is defined as
		\begin{equation}
			\mathcal{J}^{\alpha} f(t) =  \frac{1}{\Gamma(\alpha)}
			\int_a^{t} \left( \ln \frac{t}{\tau }\right)^{\alpha-1} f(\tau) \frac{\mathrm d\tau}{\tau}. 
		\end{equation}	
	\end{Definition}
	Before constructing the corresponding derivative operator we must define the following
	space of functions.
	\begin{Definition}
		Let $[a,b]$ be a finite interval such that $-\infty < a<b < \infty$ and let $AC[a,b]$ be
		the space of absolutely continuous functions on $[a,b]$. Let us denote $\delta= t\frac{\mathrm d}{\mathrm dt}$
		and define the space
		\begin{align}
    		AC_\delta^n[a,b] = \left\{ f \colon t\in [a,b] \rightarrow \mathbb{R}
    		\mbox{ \ such that \ }\left(\delta^{n-1} f\right) \in AC[a,b] \right\}.
    	\end{align}
    \end{Definition}
    Clearly $AC_\delta^1[a,b] \equiv AC[a,b]$ for $n = 1$.
     	
	\begin{Definition}
		\label{Derivative}
		Let $\delta = t\frac{\mathrm d}{\mathrm dt}$, $\Re(\alpha)>0$
		and $n = [\alpha]+1$, where $[\alpha]$ is the integer part of $\alpha$. 
		The Hadamard fractional derivative of order $\alpha$ applied to the function
		$f(t) \in AC_\delta^n[a,b]$, $0\leq a<b<\infty$, is defined as
		\begin{equation}\label{r-l}
			\mathcal{D}^{\alpha}f(t)= \frac{1}{\Gamma(n-\alpha)}\left(t\frac{\mathrm d}{\mathrm dt}\right)^n \int_a^{t}
			\left( \ln \frac{t}{\tau }\right)^{n-\alpha-1} f(\tau) \frac{\mathrm d\tau}{\tau}=\delta^n \left(\mathcal{J}^{n-\alpha}f \right) (t).
		\end{equation}	
	\end{Definition} 
    		
	It has been proved (see e.g.\ Theorem 4.8 in \cite{kilbas})
	that in the $L^p[a,b]$ space, $p\in [1,\infty)$, $0\leq a<b<\infty$, the Hadamard fractional
	derivative is the left-inverse operator to the
	Hadamard fractional integral, i.e.\
	\begin{equation}
		\label{comp}
		\mathcal{D}^{\alpha}\mathcal{J}^{\alpha}f(t)= f(t), \quad \forall t\in [a,b].
	\end{equation}
  	Analogously to the Caputo fractional calculus, the regularized Caputo-type Hadamard fractional derivative
	is defined in terms of the Hadamard fractional integral in the following way (see, for example, \cite{jara})
   	\begin{equation}
       	\left(t\frac{d}{dt}\right)^{\alpha}f(t)= \frac{1}{\Gamma(n-\alpha)}
        \int_{a}^t 
       	\left(\ln\frac{t}{\tau}\right)^{n-\alpha-1}\left(\tau \frac{d}{d\tau}\right)^n
       	f(\tau)\frac{d\tau}{\tau} = \mathcal{J}^{n-\alpha}\left(\delta^n f \right) (t),
   	\end{equation}
   	where $t\in[a,b]$, $0\leq a<b<\infty$ and $n-1<\alpha\leq n$, with $n\in \mathbb{N}$. 
	In this paper we will use the symbol 	$\displaystyle\left(t\frac{d}{dt}\right)^{\alpha}$ for the Caputo-type derivative in order to distinguish
	this one from the Riemann-Liouville type definition  \eqref{r-l} and also to underline the fact that essentially it coincides with the fractional
	power of the operator $\delta= t \ \partial/\partial t$. Moreover, by definition, when $\alpha = n$,
	$\displaystyle\left(t\frac{d}{dt}\right)^{\alpha}\equiv \delta^n$. The relationship between the Hadamard derivative \eqref{r-l} and the
	regularized Caputo-type derivative is given by (\cite{jara}, eqn.\ 12)
	\begin{equation}
    	\left(t \frac{d}{dt}\right)^\alpha f(t) = \mathcal{D}^\alpha\bigg[f(t)-\sum_{k=0}^{n-1}\frac{\delta^k f(t_0)}{k!}
    	\ln^k\left(\frac{t}{t_0}\right)\bigg], \quad \alpha \in (n-1, n], \ n = [\alpha]+1.
	\end{equation}
    In the sequel we will use the following useful equalities (that can be checked by simple calculations),
	\begin{align}
    	\left( t\frac{\mathrm d}{\mathrm dt} \right)^\alpha t^\beta = \beta^\alpha t^\beta,\\
        \mathcal{J}^\alpha t^\beta = \beta^{-\alpha} t^\beta,\label{pj}
    \end{align}
    for $\beta\in(-1,\infty)\setminus \{0\}$ and $\alpha >0$.
    It is immediate to see that 
    \begin{equation}
	    \left(t\frac{d}{dt}\right)^{\alpha}\text{\textit{const}} = 0.
    \end{equation}

%%%%%%%%%%%%%%%%%%%%%%%%%%%%%%%%%%%%%%%%%%
\section{Generalized COM-Poisson processes involving the $\alpha$-Mittag--Leffler function and the related Hadamard equations}

	Here we show a possible application of the $\alpha$-Mittag--Leffler function \eqref{0} in the context of \textit{fractional}-type Poisson
	statistics and we discuss the relation with Hadamard equations.
      	   
    We first recall that weighted Poisson distributions have been widely studied in statistics in order to consider over or under-dispersion with
    respect to the homogeneous Poisson process.
    The probability mass function of a weighted Poisson process $N^w(t)$, $t>0$, is given by (see \cite{Bala})
    \begin{equation}
		P\{N^w(t)=n\}= \frac{w(n)p(n,t)}{\mathbb{E}[w(N)]},\quad n\geq 0,
   	\end{equation}
    where $N$ is a Poisson distributed random variable 
    $$p(n,t)= \frac{t^n}{n!}e^{-t}, \quad n\geq 0, \: t>0,$$
    $w(\cdot)$ is a non-negative weight function with non-zero, finite
    expectation, i.e.
    \begin{equation}
    	0<\mathbb{E}[w(N)]= \sum_{k=0}^{\infty}w(k) p(k,t) <+\infty.
    \end{equation}
           
    The so-called COM-Poisson distribution introduced by Conway and Maxwell \cite{con} in a queueing model belongs to this wide class of
    distributions. A random variable $N_\nu(t)$ is said to have a COM-Poisson distribution if
    \begin{equation}\label{com}
	    P\{N_\nu(t)= n\}= \frac{1}{f(t)}\frac{t^n}{n!^\nu}, \quad n\geq 0, \: t>0,
    \end{equation}
    where
    \begin{equation}
    	f(t)  = \sum_{k=0}^\infty\frac{t^k}{k!^\nu}, \qquad \nu >0,
    \end{equation}
    i.e.\ the Le-Roy function.
            
	This distribution can be viewed as a particular weighted Poisson distribution with $w(n) = 1/(n!)^{\nu-1}$. The Conway-Maxwell-Poisson
	distribution is widely used in statistics in order to take into account over-dispersion (for $0<\nu<1$) or under-dispersion (for $\nu>1$)
	in data analysis. Moreover it represents a useful generalization of the Poisson distribution (that is obtained for $\nu= 1$) for many
	applications (see for example \cite{Rodado} and references therein). Observe that when $t\in (0,1)$ and $\nu \rightarrow 0$ it reduces to the
	geometric distribution, while for $\nu\rightarrow +\infty$ it converges to the Bernoulli distribution. 
			      
    Another class of interesting weighted Poisson distributions are the \textit{fractional} Poisson-type distributions
    \begin{equation}\label{mo}
		P\{N_\alpha(t)  = k  \} = \frac{1}{E_{\alpha,1}(t)}\frac{t^k}{\Gamma(\alpha k +1)}, \quad \alpha\in (0,1),
    \end{equation}
    where
    \begin{equation}
    	E_{\alpha,1}(t)= \sum_{k=0}^\infty \frac{t^k}{\Gamma(\alpha k +1)},
    \end{equation}
	is the Mittag--Leffler function (see the recent monograph \cite{gorenflo} for more details). In this case the weight function is given by
	$w(k) = k!/\Gamma(\alpha k+1)$. There is a recent wide literature about fractional generalizations of Poisson processes
	(see e.g. \cite{macci,laskin,fra,Herr,costa,bruno,scalas,fede,valla}) and we should explain in which sense we can speak about fractionality
	in the last case. The relationship of this kind of distributions with fractional calculus, first considered by Beghin and Orsingher in
	\cite{beghin}, is given by the fact that the probability generating function 
    \begin{equation}
	    G(u,t) = \sum_{k=0}^\infty u^k P\{N_\alpha(t)  = k  \}=\frac{E_{\alpha,1}(ut)}{E_{\alpha,1}(t)}, \quad \alpha \in (0,1],
    \end{equation}
    satisfies the fractional differential equation
    \begin{equation}
    	\frac{d^\alpha G(u^\alpha,t)}{du^\alpha} = t G(u^\alpha, t), \quad \alpha \in (0,1], \  |u|\leq 1,
    \end{equation}
    where $d^\alpha/du^\alpha$ denotes the Caputo fractional derivative w.r.t. the variable $u$ (see \cite{kilbas} for the definition).
    Essentially the connection between the distribution \eqref{mo} and fractional calculus is therefore given by the normalizing function,
    since the Mittag--Leffler plays the role of the exponential function in the theory of fractional differential equations.
    In this context, we discuss a new generalization of the COM-Poisson distribution by using the generalized 
    $\alpha$-Mittag--Leffler function \eqref{0}, independently introduced by Gerhold in \cite{ge} and by Garra and Polito in \cite{had}.
	This new \textit{fractional} distribution includes a wide class of special distributions already studied in the recent literature:
    the fractional Poisson distribution \eqref{mo} considered by Beghin and Orsingher \cite{beghin}, the classical 
    COM-Poisson distribution and so on.
   
    In this paper we introduce the following generalization, namely the \textit{fractional} COM-Poisson distribution
    \begin{equation}\label{a}
    	P\{N_{\nu;\alpha,\gamma}(t)=n\}=\frac{(\lambda t)^n}{\left(\Gamma(\alpha n+\gamma)\right)^\nu}\frac{1}{E_{\nu; \alpha,\gamma}(\lambda t)},
    \end{equation}
    where 
    \begin{equation}
    	E_{\nu; \alpha,\gamma}(t)= \sum_{k=0}^{\infty}
        \frac{t^k}{\left(\Gamma(\alpha k+\gamma)\right)^\nu}, \quad \alpha >0, \ \nu>0, \ \gamma \in \mathbb{R}
    \end{equation}
    is the $\alpha$-Mittag--Leffler function, considered by Gerhold  \cite{ge}, Garra and Polito \cite{had}. The distribution \eqref{a} is a weighted
    Poisson distribution with weights $w(k)= k!/\left(\Gamma(\alpha k + \gamma)\right)^\nu$.
    Obviously the probability generating function is here given by
    a ratio of $\alpha$-Mittag--Leffler functions, i.e. 
    \begin{equation}
	    G(u,t)= \sum_{n=0}^{\infty}u^n  P\{N_{\nu;\alpha,\gamma}(t)=n\}=
        \frac{E_{\nu; \alpha,\gamma}(\lambda u  t)}{E_{\nu; \alpha,\gamma}(\lambda t)}, \quad |u| \leq 1.
    \end{equation}
    With the next theorem we explain the reason why 
    we consider the distribution \eqref{a}, a kind of \textit{fractional}
    generalization of the COM-Poisson classical distribution related to Hadamard fractional equations with varying coefficients.
          
    \begin{Theorem}
	    For $\nu \in(0,1)$, the function
        $$g(u,t)=u^{\nu-1}\frac{E_{\nu;1, \nu}(\lambda ut)}{E_{\nu;1, \nu}(\lambda t)}= u^{\nu-1} G(u,t),$$
        satisfies the equation
        \begin{equation}
    	    \left(u\frac{d}{du}\right)^{\nu}g(u,t)=
            \lambda t u g(u,t)+\frac{u^{\nu-1}}{\Gamma^\nu(\nu-1)} ,
        \end{equation}
        where $\displaystyle{\left(u\frac{d}{du}\right)^{\nu}}$ is the 
        Caputo-Hadamard fractional derivative of order $\nu$.
    \end{Theorem}
    \begin{proof}
	    By direct calculations we obtain that
    	\begin{align}
	      	\left(u\frac{d}{du}\right)^{\nu}u^{\nu-1}\frac{E_{\nu;1,\nu}(\lambda u t)}{E_{\nu;1,\nu}(\lambda  t)}
	      	&= \frac{1}{E_{\nu;1,\nu}(\lambda  t)}\sum_{k=0}^{\infty}
	   		\frac{\lambda^k u^{k+\nu-1}t^k}{\left(\Gamma(k+\nu-1)\right)^\nu}\\
	   		\nonumber& = \frac{1}{E_{\nu;1,\nu}(\lambda  t)} 
	   		\sum_{k=-1}^{\infty}
	 		\frac{\lambda^{k+1} u^{k+\nu}t^{k+1}}{\left(\Gamma(k+\nu)\right)^{\nu}}
	  		= \lambda t u g(u)+\frac{u^{\nu-1}}{\left(\Gamma(\nu-1)\right)^{\nu}},
		\end{align}   
   		where we changed the order of summation with fractional differentiation since the $\alpha$-Mittag--Leffler function
    	is, for $\alpha \in (0,1)$, an entire function.
    \end{proof}
        	          
    A more detailed analysis about the statistical properties and possible applications of this new class
    of distributions is not object of this paper and will be considered in a future work.
             
    We finally observe that the fractional COM-Poisson distribution \eqref{a} includes as a special case, for 
    $\nu = 1$ the fractional Poisson distribution \eqref{mo} previously introduced by Beghin and Orsingher in \cite{beghin}
    and recently treated by Chakraborty and Ong \cite{cha},  Herrmann \cite{Herr} and Porwall and Dixit \cite{dix}. 
    Asymptotic results for the multivariate version of this distribution have been recently analyzed by Beghin and Macci in \cite{macci}.

%%%%%%%%%%%%%%%%%%%%%%%%%%%%%%%%%%%%%%%%%%
\section{Hadamard fractional equations related to the GCOM-Poisson distribution}

	In a recent paper, Imoto \cite{imoto} studied a different generalization of the Conway-Maxwell-Poisson distribution (in the following GCOM Poisson
	distribution), depending from two real parameters $r$ and $\nu$. According to \cite{imoto}, a random variable $X(t)$ is said to have a GCOM Poisson
	distribution if
    \begin{equation}\label{gc}
	    P\{X(t)= n\}= \frac{\left(\Gamma(\nu+n)\right)^r t^n}{n!C(r,\nu,t)}, \quad t\geq 0, \ n= 0,1, \dots,
    \end{equation}
    where 
    \begin{equation}\label{cno}
    	C(r,\nu, t)= \sum_{k=0}^\infty \frac{\left(\Gamma(\nu+k)\right)^r t^k}{k!},
    \end{equation}
    is the normalizing function which converges for $\nu >0$, $t>0$ and $r<1/2$ (by applying the ratio criterium) or $r=1$, $\nu>0$ and $|t|<1$.
    In the context of weighted Poisson distributions, it corresponds to the choice $w(k) = \displaystyle \Gamma^r(\nu+k)$.
    
    This distribution includes as particular cases the COM-Poisson distribution for $\nu = 1$ and $r = 1-n$, the geometric for $r=\nu = 1$, $|t|<1$
    and the homogeneous Poisson distribution for $r=0$. The new parameter introduced in the distribution plays the role
    of controlling the length of tails within the framework of queueing processes. Indeed it was proved that, it gives over-dispersion for
    $r \in(0,1)$ and under-dispersion for $r<0$. Moreover, from the ratios of successive probabilities, when $\nu >1$ and $r\in(0,1)$ it is heavy
    tailed, while for $r<0$ light-tailed.
          
    We also underline that for $r = 1-n$, with $n>0$, the function \eqref{cno} becomes a sort of generalization of the Wright function
    in the spirit of the $\alpha$-Mittag--Leffler function
    \begin{equation}\label{cno1}
	    C(1-n,\nu, t)= \sum_{k=0}^\infty  \frac{t^k}{k!\Gamma^n(\nu+k)}.
    \end{equation}
    Recall that, in the general form, the Wright function is defined as follows (we refer to \cite{gim} for a good survey):
    \begin{equation}
        W_{\alpha,\beta}(t) = \sum_{k=0}^\infty \frac{t^k}{k!\Gamma(\alpha k +\beta)}, \quad \alpha >-1, \beta\in \mathbb{C} .
    \end{equation}
         
    We show now the relation between the normalizing function \eqref{cno} and an Hadamard fractional differential equation with varying coefficients.
    Therefore, we can consider in some sense also the GCOM Poisson distribution as a \textit{fractional}-type modification of the homogeneous
    distribution, in the sense that its probability generating function satisfies a fractional differential equation.    
     
    \begin{Proposition}
    	The function 
     	\begin{equation}
     		C(r,\nu,t)= \sum_{k=0}^\infty \frac{\Gamma^r(\nu+k) (\lambda t)^k}{k!},
     	\end{equation}
     	satisfies the  integro-differential equation 
     	\begin{equation}\label{ama}
     		\mathcal{J}^{r}\left(t^\nu \frac{df}{dt}\right)= \lambda t^\nu f, \quad t>0, \ \nu\in(0,1], \ r\in (0,1/2)
     	\end{equation}
     	involving an Hadamard fractional integral $\mathcal{J}^{r}$ of order $r$.
    \end{Proposition}
     
    \begin{proof}
    	\begin{align}
      		\mathcal{J}^{r}\left(t^\nu \frac{dC}{dt}\right)&=\mathcal{J}^{r} \sum_{k=1}^\infty \frac{\Gamma^r(\nu+k) \lambda^k t^{k+\nu-1}}{(k-1)!}= 
      		\sum_{k=1}^\infty \frac{\Gamma^r(\nu+k) (k+\nu-1)^{-r} \lambda^k t^{k+\nu-1}}{(k-1)!}\\
      		\nonumber &  = \sum_{k=1}^\infty \frac{\Gamma^r(\nu+k-1) \lambda^k t^{k+\nu-1}}{(k-1)!} = \lambda t^\nu f,
     	\end{align}
     	where we used \eqref{pj}.
    \end{proof}
     
    \begin{Corollary}
     	The probability generating function $\mathcal{G}(u,t)$ of the GCOM Poisson distribution \eqref{gc} 
    	\begin{equation}
     		\mathcal{G}(u,t)= \sum_{k=0}^\infty u^k P\{X(t)=k\} =  
     		\frac{C(r,\nu, u t)}{C(r,\nu,t)}, \quad |u|\leq 1,
     	\end{equation}
     	satisfies the integro-differential equation
     	\begin{equation}
     		\mathcal{J}^{r}\left(u^\nu \frac{df}{du}\right)= u^\nu f,\quad |u|\leq 1. 
     	\end{equation}
    \end{Corollary}

\section{Further connections between modified Mittag--Leffler functions and Hadamard fractional equations.}
   
	The analysis of fractional differential equations with varying coefficients is an interesting and non-trivial topic. In particular the analysis of
	equations involving fractional-type Bessel operators (i.e. the fractional counterpart of singular linear differential operators of arbitrary order)
	has attracted the interest of many researchers (see e.g. \cite{virginia,sitnik} and references therein). In \cite{had}, some results about the
	connection between Le-Roy functions and equations with space-varying coefficients involving Hadamard derivatives and Laguerre derivatives have
	been obtained. 
	Here we go on in the direction started in \cite{had}, showing other interesting applications of Hadamard fractional equations in the theory of
	hyper-Bessel functions. With the next theorem, we find the equation interpolating classical Bessel equations of arbitrary order possessing exact
	solution in terms of Le Roy functions. We remark, indeed, that Le Roy functions include as special cases hyper-Bessel functions.
   
    \begin{Theorem}
	    The Le Roy function $E_{\alpha;1,1} \left(t^\alpha/\alpha \right)$ satisfies the equation
        \begin{equation}\label{bilb}
        	\frac{1}{t^{\alpha}}\left(t\frac{d}{dt}\right)^{\alpha}f(x)= \alpha^{\alpha-1}f(t)
        \end{equation}
	\end{Theorem}
  		 
    \begin{proof}
    	Recall that
        \begin{align}
        	\left( t \frac{\mathrm d}{\mathrm dt} \right)^\alpha t^\beta = \beta^\alpha t^\beta.
        \end{align}
        Therefore
        \begin{align}
        	\left( t \frac{\mathrm d}{\mathrm dt} \right)^\alpha
            E_{\alpha;1,1} \left( \frac{t^\alpha}{\alpha} \right)
            = \sum_{k=0}^\infty \frac{(\alpha k)^\alpha t^{\alpha k}}{\alpha^k k!^\alpha}
            = \alpha^{\alpha-1} t^\alpha E_{\alpha;1,1} \left( \frac{t^\alpha}{\alpha} \right).
        \end{align}
	\end{proof}
    
    \begin{Remark}
	    For $\alpha=2$, we have that the function $E_{2;1,1}(t^2/2)$ is a solution
        of
        \begin{align}
    	    \left( t\frac{\mathrm d}{\mathrm dt} \right)^2 f(t) = 2 t^2 f(t).
        \end{align}
        For $\alpha=3$, the function $E_{3;1,1}(t^3/3)$ satisfies
        \begin{align}
        	\left( t\frac{\mathrm d}{\mathrm dt} \right)^3 f(t) = 9 t^3 f(t),
        \end{align}
        and so forth for any integer value of $\alpha$. From this point of view, the solution of Hadamard equation \eqref{bilb}, for non integer
        values of $\alpha$, leads to an interpolation between successive hyper-Bessel functions.
	\end{Remark}
       
    Let us consider the operator
    \begin{align}
	    L_H = \frac{1}{t^{\alpha n}} \underbrace{\left( t \frac{d}{dt} \right)^\alpha \left( t \frac{d}{dt} \right)^\alpha
	    \dots \biggl( t \frac{d}{dt} \biggr)^\alpha}_{\text{$n$ times}}, \quad \alpha >0,
    \end{align}
    where $H$ stands for an Hyper-Bessel type operator involving Caputo-type Hadamard derivatives. We have the following
    general
    \begin{Theorem}
    	A solution of the equation
    	\begin{equation}
    		L_H f(t)= \alpha^{n\alpha-n} n^{n\alpha}f(t),
    	\end{equation}
       is given by
       $$f(t)= E_{n\alpha;1,1} \left( t^{\alpha n}/\alpha^n \right)$$
	\end{Theorem}
  
   	To conclude this section, we restate Theorem 3.3 of \cite{had}, in view of the comments of Turmetov \cite{turmetov} in the case in which
   	Caputo-type Hadamard derivatives appear in the governing equations.
	\begin{Theorem}
		The function $E_{\beta;1,1}(t)$, with $\beta\in [1,\infty)$, $t\geq 0 $, $\lambda \in \mathbb{R}$ is an eigenfunction of the operator
        \begin{align}
        	\label{caputoite}
            \underbrace{\frac{\mathrm d}{\mathrm dt} t \frac{\mathrm d}{\mathrm dt}
            \dots \frac{\mathrm d}{\mathrm dt} t
            \frac{\mathrm d}{\mathrm dt}}_{\text{$r$ derivatives}}\left(t\frac{d}{dt}\right)^{\beta-r}, \quad r = 1,\dots, n-1, 
        \end{align}
        where $n = [\beta]$ is the integer part of $\beta$ and $\displaystyle \left(t\frac{d}{dt}\right)^{\beta-r}$ denotes the Caputo-type
        regularized Hadamard derivative of order $\beta-r$.
    \end{Theorem}
   
    The difference w.r.t.\ the previous version is simply given by the fact that, by using the regularized Caputo-type Hadamard
    derivative $\displaystyle\left(t\frac{d}{dt}\right)^{\alpha}t^0 = 0$, that is necessary for the correctness of the result.

%%%%%%%%%%%%%%%%%%%%%%%%%%%%%%%%%%%%%%%%%%
\acknowledgments{The work of R.G. has been carried out in the framework of the activities of GNFM. F.P.\ has been supported
		by the projects \emph{Memory in Evolving Graphs} (Compagnia di San Paolo/Universit\`a di Torino),
		\emph{Sviluppo e analisi di processi Markoviani e non Markoviani con applicazioni} (Universit\`a di Torino), and by INDAM--GNAMPA.}

%%%%%%%%%%%%%%%%%%%%%%%%%%%%%%%%%%%%%%%%%%
%\authorcontributions{For research articles with several authors, a short paragraph specifying their individual contributions must be provided. The following statements should be used ``X.X. and Y.Y. conceived and designed the experiments; X.X. performed the experiments; X.X. and Y.Y. analyzed the data; W.W. contributed reagents/materials/analysis tools; Y.Y. wrote the paper.'' Authorship must be limited to those who have contributed substantially to the work reported.}

%%%%%%%%%%%%%%%%%%%%%%%%%%%%%%%%%%%%%%%%%%
%\conflictsofinterest{The authors declare no conflict of interest.} 

%%%%%%%%%%%%%%%%%%%%%%%%%%%%%%%%%%%%%%%%%%
%% optional
%\abbreviations{The following abbreviations are used in this manuscript:\\

%\noindent 
%\begin{tabular}{@{}ll}
%MDPI & Multidisciplinary Digital Publishing Institute\\
%DOAJ & Directory of open access journals\\
%TLA & Three letter acronym\\
%LD & linear dichroism
%\end{tabular}}

%%%%%%%%%%%%%%%%%%%%%%%%%%%%%%%%%%%%%%%%%%
%% optional


%=====================================
% References, variant A: internal bibliography
%=====================================
\reftitle{References}

% The following MDPI journals use author-date citation: Arts, Econometrics, Economies, Genealogy, Humanities, IJFS, JRFM, Laws, Religions, Risks, Social Sciences. For those journals, please follow the formatting guidelines on http://www.mdpi.com/authors/references
% To cite two works by the same author: \citeauthor{ref-journal-1a} (\citeyear{ref-journal-1a}, \citeyear{ref-journal-1b}). This produces: Whittaker (1967, 1975)
% To cite two works by the same author with specific pages: \citeauthor{ref-journal-3a} (\citeyear{ref-journal-3a}, p. 328; \citeyear{ref-journal-3b}, p.475). This produces: Wong (1999, p. 328; 2000, p. 475)

%=====================================
% References, variant B: external bibliography
%=====================================
%\externalbibliography{yes}
%\bibliography{your_external_BibTeX_file}

%%%%%%%%%%%%%%%%%%%%%%%%%%%%%%%%%%%%%%%%%%
%% optional
%\sampleavailability{Samples of the compounds ...... are available from the authors.}

%%%%%%%%%%%%%%%%%%%%%%%%%%%%%%%%%%%%%%%%%%
\end{document}